\documentclass[mnsc,nonblindrev]{informs3_hide} 

\OneAndAHalfSpacedXI 

\usepackage{amsmath,amsfonts,amssymb}
\usepackage{mathtools}
    \allowdisplaybreaks
\usepackage{bm}
\usepackage[mathscr]{euscript}
\usepackage{enumitem}
\usepackage{xcolor}
\usepackage{graphicx}
\usepackage{multirow}
\usepackage{booktabs}
\usepackage{microtype}
\usepackage{fix-cm}
\usepackage[ruled,vlined]{algorithm2e}
\usepackage{physics}
\usepackage{xfrac}
\usepackage{nicematrix}

\DeclareMathOperator{\E}{\mathbb E}
\DeclareMathOperator{\Var}{\mathrm{Var}}
\DeclareMathOperator{\Cov}{\mathrm{Cov}}

\DeclareMathOperator{\Real}{\mathbb R}

\newcommand{\iu}{\mathrm{i}\mkern1mu}

\usepackage{natbib}
 \bibpunct[, ]{(}{)}{,}{a}{}{,}%
 \def\bibfont{\small}%
\TheoremsNumberedThrough     
\ECRepeatTheorems

\EquationsNumberedThrough    

\MANUSCRIPTNO{}

\begin{document}


\RUNAUTHOR{Hong and Zhang}

\RUNTITLE{Surrogate-Based Simulation Optimization}

\TITLE{Surrogate-Based Simulation Optimization}

\ARTICLEAUTHORS{
\AUTHOR{L. Jeff Hong}
\AFF{School of Management and School of Data Science, Fudan University, Shanghai 200433, China, \EMAIL{hong\_liu@fudan.edu.cn}}
\AUTHOR{Xiaowei Zhang}
\AFF{Faculty of Business and Economics, The University of Hong Kong, Pokfulam Road, Hong Kong SAR, \EMAIL{xiaoweiz@hku.hk}}
}

\ABSTRACT{Simulation models are widely used in practice to facilitate decision-making in a complex, dynamic and stochastic environment.
But they are computationally expensive to execute and optimize, due to lack of analytical tractability.
Simulation optimization is concerned with developing efficient sampling schemes---subject to a computational budget---to solve such optimization problems.
To mitigate the computational burden, surrogates are often constructed using simulation outputs to approximate the response surface of the simulation model.
In this tutorial, we provide an up-to-date overview of surrogate-based methods for simulation optimization with continuous decision variables.
Typical surrogates, including linear basis function models and Gaussian processes, are introduced.
Surrogates can be used either as a local approximation or a global approximation.
Depending on the choice, one may develop algorithms that converge to either a local optimum or a global optimum.
Representative examples are presented for each category.
Recent advances in large-scale computation for Gaussian processes are also discussed.
}

\KEYWORDS{Surrogate, simulation optimization, Gaussian process, matrix inversion}


\maketitle

\section{Introduction}

Simulation optimization (SO) concerns a class of optimization problems whose objective functions and/or constraints do not possess an analytical form and can only be evaluated based
on noisy simulation samples.
The simulation model is usually expensive to execute;
thus, the number of evaluations that one is allowed to perform is limited, subject to one's computational budget.
Specifically, we consider in this tutorial problems of the form
\begin{equation}\label{eq:SO}
\max_{\BFx\in\mathscr{X}} \{f(\BFx) \coloneqq \E[F(\BFx)]\},
\end{equation}
where $\BFx$ is the decision variable, $\mathscr{X}$ is the feasible set, and $F(\BFx)$ is a real-valued random variable, representing the stochastic response of a simulation model evaluated at $\BFx$.
The distribution of $F(\BFx)$ is an unknown function of $\BFx$,
but its samples can be generated from running simulation experiments.

Depending on the nature of the feasible set, problem~\eqref{eq:SO} demands a distinct treatment and principle for designing algorithms.
When $\mathscr{X}$ is a set of a relatively small number of feasible solutions with no inherent ordering defined,
the problem is known as \emph{ranking and selection} (R\&S).
For example, $\mathscr{X}$ may represent feasible system configurations  with regard to how many and what redundant components to use to design a reliable system.

Another common setting is that $\mathscr{X}$ is integer-ordered, that is, $\mathscr{X} = \Omega \cap \mathbb{Z}^d$, where $\Omega \subset \Real^d$ is a convex set and
$\mathbb{Z}^d$ is the set of $d$-dimensional integer vectors.
In this setting, $\mathscr{X}$ usually has a very large or even  infinite number of elements,
and problem~\eqref{eq:SO} is called \emph{discrete optimization via simulation} (DOvS).
For example, a retailer may need to make stocking decisions for $d$ products to minimize operational costs,  $\BFx$ may represent the numbers of units to order for these products,
and $\mathscr{X}$ may be formed by capacity constraints.

In this tutorial, we assume that $\BFx$ is a $d$-dimensional real vector and $\mathscr{X}\subseteq \Real^d$ is a continuous set,
and refer to  \cite{HongNelsonXu15} for a recent tutorial on R\&S and DOvS problems.
However,
our concentration on continuous decision variables does not limit the scope of this tutorial.
In general, algorithms designed for
the continuous setting can  be applied to DOvS, although theoretical analyses may need to be adjusted.

A variety of strategies  have been proposed to solve continuous SO problems,
including sample average approximation \citep{KimPasupathyHenderson15},
stochastic approximation \citep{ChauFu15},
and random search \citep{Andradottir15}.
This tutorial focuses on surrogate-based methods---another class of strategies---which have gained increasing interest in recent years, despite their long history.
The popular use of surrogates in SO may be attributed to (i) the flexibility to capture complex surface shapes and (ii) the capability to predict surface values where no simulation samples are observed.
The latter reason is of particular importance in light of the fact that computational budget is generally regarded as a scarce resource relative to the cost of simulation experiments.

Finally, we note that in contrast to earlier articles on surrogate-based methods \citep{BartonMeckesheimer06,Barton09},
the present tutorial aims to bring the readers up to date on the fast developments in the area.
One featured discussion is on recent advances in coping with the computational challenges associated with the use of surrogates for large-scale problems.


The rest of this tutorial is organized as follows.
Section~\ref{sec:surrogate} introduces surrogates that are widely used in practice.
Section~\ref{sec:local} and Section~\ref{sec:global} present SO methods that use surrogates as local approximations and global approximations, respectively.
Section~\ref{sec:largescale} discusses recent advances in handling computational issues that arise when using surrogates for large datasets.
Section~\ref{sec:conclusions} highlights current research challenges and potential opportunities.

\section{Surrogates}\label{sec:surrogate}

A surrogate---also known as (a.k.a.) metamodel---is an approximation to the response surface, that is, the simulation input-output relationship.
The main purpose of using a surrogate is to mitigate the computational burden of running expensive simulation experiments.
Although any supervised learning \citep{HastieTibshiraniFriedman09} model may, in principle, be used,
one typically has several considerations to keep in mind  when choosing a surrogate to cope with computational budget constraints.
\begin{enumerate}[label=(\roman*)]
\item It should possess a simple structure and does not require ``big data'' to fit, because simulation samples are expensive to acquire.
\item It should be computationally easy to fit, because it often needs to be updated in a sequential fashion as more simulation samples become available.
\item It should give rise to a predictor in explicit form, so that
predictions can be computed efficiently, theoretical analysis can be facilitated, and the surrogate can be optimized easily.
\end{enumerate}

Three classes of surrogates that satisfy the above criteria have been widely adopted in practice: low-order polynomials, linear basis function models, and Gaussian processes; we introduce them in Sections~\ref{sec:polynomial}, \ref{sec:basis}, and \ref{sec:GP}, respectively.
The first---including linear and quadratic functions---are normally viewed as local approximations from the perspective of Taylor expansion,
whereas the last two are global approximations thanks to their nonparametric nature.
Nevertheless, as we demonstrate in Section~\ref{sec:connect}, they can be unified through the lens of ridge regularization.
In Section~\ref{sec:enhance}, we present recent approaches to enhancing the prediction capability of a surrogate if additional information is available.

Before formally presenting the surrogates, we state the common set-up.
Suppose that the simulation model is executed at $\{\BFx_i:i=1,\ldots,n\}$,
where $\BFx_i=(x_{i,1},\ldots,x_{i,d})$ for each $i$.
For each $\BFx_i$, the number of replications is $r_i\geq 1$; each replication generates a realization of the random variable $F(\BFx_i)$ in Equation~\eqref{eq:SO}, denoted by  $y_{i,\ell}$ for $\ell=1,\ldots,r_i$.
Let $y_{i,\ell} = f(\BFx_i) + \epsilon_{i,\ell}$,
where $\epsilon_{i,\ell}$ is independent Gaussian noise with mean zero and variance $\sigma^2(\BFx_i)$.
We are interested in approximating $f$ via the regression equation
\begin{equation}\label{eq:generic_reg}
\bar{y}_i = f(\BFx_i) + \bar{\epsilon}_i, \quad i=1,\ldots,n,
\end{equation}
where $\bar{y}_i \coloneqq {r_i}^{-1}\sum_{\ell=1}^{r_i} y_{i,\ell}$
and $\bar{\epsilon}_i \coloneqq {r_i}^{-1}\sum_{\ell=1}^{r_i} \epsilon_{i,\ell} $.
Clearly, $\Var[\bar{\epsilon}_i] = \sigma^2(\BFx_i)/r_i$.

\subsection{Polynomials}\label{sec:polynomial}

Due to the explosion in the number of terms in the representation of a polynomial in multiple dimensions,
polynomials with orders higher than two are seldomly used to approximate a response surface.
Low-order polynomials are suitable for situations where we are interested in a localized region of the feasible set (or design space) $\mathscr{X}$.
A second-order polynomial (i.e., quadratic function) is
\begin{equation}\label{eq:quadratic}
f(\BFx) = \beta_0 + \sum_{j=1}^d \beta_j x_j + \sum_{j=1}^d\sum_{k=1}^d \beta_{jk}x_jx_k,
\end{equation}
where $\BFx = (x_1,\ldots,x_d)$.
This surrogate may be appropriate if we expect the response surface $f$ has substantial curvature.
In contrast, a first-order polynomial (i.e., linear function), to which Equation~\eqref{eq:quadratic} is reduced by setting $\beta_{jk}=0$ for all $j$ and $k$, may be a better fit in the presence of little curvature.

The parameters $\beta_0, \beta_j, \beta_{jk}$ can be estimated via ordinary least squares (OLS)
\begin{equation}\label{eq:OLS}
\min_{\beta_0,\beta_j,\beta_{jk}} \frac{1}{n} \sum_{i=1}^n \biggl(\bar{y}_i - \beta_0 - \sum_{j=1}^d \beta_j x_{i,j} - \sum_{j=1}^d\sum_{k=1}^d \beta_{jk}x_{i,j}x_{i,k} \biggr)^2.
\end{equation}
The solution is given in Section~\ref{sec:basis} in a more general setting.
Then, one can predict the response at any arbitrary location $\BFx$
by plugging the OLS estimates $\hat\beta_0, \hat\beta_j, \hat\beta_{jk}$
in Equation~\eqref{eq:quadratic}, that is,
\[
\hat f(\BFx) = \hat\beta_0 + \sum_{j=1}^d \hat\beta_j x_{j} + \sum_{j=1}^d\sum_{k=1}^d \hat \beta_{jk}x_{j}x_{k}.
\]

\subsection{Linear Basis Function Models}\label{sec:basis}

If low-order polynomials do not provide a good fit,
a natural extension uses linear basis function models that express
$f$ as a linear combination of basis functions,
\begin{equation}\label{eq:lienar-basis}
f(\BFx) = \BFbeta^\intercal \BFphi(\BFx) =  \sum_{k=1}^p \beta_k \phi_k(\BFx),
\end{equation}
where $\BFbeta = (\beta_1,\ldots,\beta_p)^\intercal$ is a vector of unknown parameters and $\BFphi(\BFx) = (\phi_1(\BFx),\ldots,\phi_p(\BFx))^\intercal$ is a vector of chosen basis functions, such as the truncated power basis and the Fourier basis, etc.; see \citet[Chapter 3]{Bishop06} and \citet[Chapter 5]{HastieTibshiraniFriedman09}.

A particular popular class, among others, is radial basis functions (RBFs), meaning that the value of $\phi$ depends on $\BFx$ only through the distance between $\BFx$ and some fixed point, say $\BFc\in\Real^d$.
That is, an RBF has the form $\varphi(\|\BFx-\BFc\|)$ for some function $\varphi:\Real\mapsto\Real$.
Typical examples include
the Gaussian RBF $\varphi(x) = \exp(-x^2/2\eta^2)$ and the thin plate spline $\varphi(x) = x^2 \ln(x)$.

The parameters $\BFbeta$ can be also estimated via OLS:
\begin{equation}\label{eq:OLS-sol}
\begin{aligned}
\hat\BFbeta ={}& \argmin_{\BFbeta} \frac{1}{n}\sum_{i=1}^n (\bar{y}_i - \BFbeta^\intercal \BFphi(\BFx_i))^2 \\
={}& \BFPhi^\intercal (\BFPhi \BFPhi^\intercal )^{-1}\bar{\BFy},
\end{aligned}
\end{equation}
where $\BFPhi$ is the $n$-by-$p$ matrix with the $i$-th row being $\BFphi(\BFx_i)^\intercal$ for all $i=1,\ldots,n$. The prediction is then given by
\begin{equation}\label{eq:basis-pred}
\hat f(\BFx) = \hat\BFbeta^\intercal \BFphi(\BFx) =  (\BFPhi \BFphi(\BFx) )^\intercal (\BFPhi \BFPhi^\intercal )^{-1}\bar{\BFy}.
\end{equation}

\subsection{Gaussian Processes}\label{sec:GP}

Gaussian processes (GPs, a.k.a. Gaussian random fields) have a remarkable success in numerous research areas.
Using them as surrogates originated in geostatistics, where the method was named \emph{kriging} \citep{Krige51,Matheron63}.
Later, they were applied to the design and analysis of (deterministic) computer experiments \citep{SacksWelchMichellWynn89}.
The adoption of GPs to approximate response surfaces in stochastic simulation literature, where the method is often called \emph{stochastic kriging}, was popularized by \cite{AnkenmanNelsonStaum10}.
The same method was referred to as \emph{GP regression} in the machine learning literature \citep{RasmussenWilliams06}.

Whereas linear regression with basis functions is a frequentist approach,
GPs represent a Bayesian viewpoint.
A GP with domain $\mathscr{X}$ is fully characterized by its mean function $\mu:\mathscr{X}\mapsto \Real$ and covariance function $K:\mathscr{X}\times \mathscr{X}\mapsto \Real$.
To use a GP as a surrogate,
we start by imposing a GP prior on the unknown function $f$, denoted by $f\sim \mathsf{GP}(\mu,K)$.
Some may prefer to state the assumption in a different way---with a somewhat less Bayesian flavor---namely, $f$ is a sample path (i.e., realization) of the GP.
Both mean the following:
for any finite set $\{\BFx_1,\ldots,\BFx_n\}\subset \mathscr{X}$ of any size $n\geq 1$,
$(f(\BFx_1),\ldots,f(\BFx_n))$ has a multivariate normal distribution with mean vector $\BFmu = (\mu(\BFx_1),\ldots,\mu(\BFx_n))^\intercal$ and $n$-by-$n$ covariance matrix $\BFK = (K(\BFx_i, \BFx_{i'}))_{i,i'=1}^n$.

\subsubsection{Mean Functions}

From a modeling point of view, the mean function $\mu$ is used to encode
one's prior knowledge about the overall shape of the response surface $f$.
It is usually chosen in one of the following ways.
\begin{enumerate}[label=(\roman*)]
\item Set $\mu(\BFx)\equiv c$ for some constant $c$ representing the overall surface mean. This is possibly the most common treatment in practice.
One may even set $c=0$ if little prior knowledge is available.
\item Set $\mu(\BFx) = \BFbeta^\intercal \BFphi(\BFx)$, where $\BFphi(\BFx)$ is a vector of known basis functions and $\BFbeta$ is a vector of  hyperparameters of compatible dimension.
\item Set $\mu(\BFx)$ to be a function derived from a simplified, analytical model of the same stochastic system that the original simulation model aims to describe.
For example, if the original simulation model is a complex queueing model,
$\mu$ may be derived in closed form based on a highly stylized queueing model; see Section~\ref{sec:enhance} for details.
\end{enumerate}

\subsubsection{Covariance Functions}\label{sec:cov}

We first introduce two popular classes of covariance functions: the Gaussian class and the Mat\'ern class.
They are both \emph{stationary}, meaning that $\Cov[f(\BFx), f(\BFx')]$ depends on $\BFx$ and $\BFx'$ only through the difference $(\BFx-\BFx')$.
We then present a class of covariance functions that permits a different  level of differentiability in each dimension.
More examples of covariance functions can be found in \citet[Chapter 4]{RasmussenWilliams06}.

\begin{example}[Gaussian Covariance Functions]
For positive constants $\tau$ and $\eta$, a Gaussian covariance function is defined by
\begin{equation}\label{eq:Gaussian-kernel}
K_{\mathsf{Gaussian}}(\BFx,\BFx') = \tau^2 \exp(-\frac{\|\BFx-\BFx'\|^2}{2\eta^2}), \quad \BFx,\BFx'\in\Real^d.
\end{equation}
It is also called a squared exponential covariance function.
\end{example}

\begin{example}[Mat\'ern Covariance Functions]
For positive constants $\tau$, $\eta$, and $\nu$, a Mat\'ern covariance function is defined by
\begin{equation}\label{eq:Matern-kernel}
K_{\mathsf{Matern}}(\BFx,\BFx';\nu) =  \frac{\tau^2 }{\Gamma(\nu)2^{\nu-1}}\left(\frac{\sqrt{2\nu} \| \BFx-\BFx'\|}{\eta}\right)^{\nu} \mathsf{K}_{\nu}\left(\frac{\sqrt{2\nu} \| \BFx-\BFx'\|}{\eta}\right), \quad \BFx,\BFx'\in\Real^d,
\end{equation}
where $\Gamma(\cdot)$ is the gamma function, and $\mathsf{K}_\nu(\cdot)$ is the modified Bessel function of the second kind of order $\nu$.
The parameter $\nu$ is usually set to be half-integer, i.e., $\sfrac{1}{2}, \sfrac{3}{2},\sfrac{5}{2},\ldots$, in which case the expression of  $K_{\mathsf{Matern}}(\BFx,\BFx';\nu)$ can be simplified substantially. For instance,
\[
K_{\mathsf{Matern}}(\BFx,\BFx';\nu) =
\left\{
\begin{array}{ll}
\displaystyle \tau^2\exp\biggl(\frac{- \|\BFx-\BFx'\|}{\eta}\biggr), &\quad\mbox{if }\nu=1/2,  \\[1ex]
\displaystyle  \tau^2\biggl(1+ \frac{\sqrt{3}\|\BFx-\BFx'\|}{\eta}\biggr)\exp\biggl(-\frac{ \sqrt{3}\| \BFx-\BFx'\|}{\eta}\biggr), &\quad\mbox{if }\nu=3/2,  \\[1ex]
\displaystyle  \tau^2\biggl(1+ \frac{\sqrt{5}\|\BFx-\BFx'\|}{\eta} + \frac{5\|\BFx-\BFx'\|^2}{3\eta^2}\biggr)\exp\biggl(-\frac{ \sqrt{5}\| \BFx-\BFx'\|}{\eta}\biggr), &\quad\mbox{if }\nu=5/2.
\end{array}
\right.
\]
A general formula can be found on page 85 of \cite{RasmussenWilliams06}.
\end{example}

\smallskip

Note that both the expressions in \eqref{eq:Gaussian-kernel} and \eqref{eq:Matern-kernel}
are in the form of $\tau^2 \rho(\|\BFx-\BFx'\|)$ for some function $\rho:\Real_{\geq 0}\mapsto(0,1]$.
Thus, the parameter $\tau^2$ represents the marginal variance of the associated GP and $\rho(\|\BFx-\BFx'\|)$ represents the correlation.

The parameter $\nu$ of the Mat\'ern covariance functions is called the \emph{smoothness} parameter,
for it controls the order of differentiability of the sample paths of the induced GP; see \citet[Section 6.5]{Stein99}.
It can be shown \citep[page 50]{Stein99} that
\[
\lim_{\nu\to\infty} K_{\mathsf{Matern}}(\BFx,\BFx';\nu) = K_{\mathsf{Gaussian}}(\BFx,\BFx'), \quad \BFx,\BFx'\in\Real^d.
\]
This suggests that the sample paths induced by the Gaussian covariance functions are infinitely differentiable---an overly strong property that may not be reasonable for some response surfaces.

Being controlled by a single parameter $\nu$, the differentiability of the GPs associated with the Mat\'ern covariance functions are homogeneous in each dimension.
Motivated by the need for flexibility to control differentiability separately in different dimensions,  \cite{SalemiStaumNelson19} propose a new class of covariance functions, and the corresponding GPs are called \emph{generalized integrated Brownian fields} (GIBFs).
Unlike the Gaussian and Mat\'ern classes, the GIBF covariance functions are nonstationary and possess a tensor product form.
\begin{example}[GIBF Covariance Functions]
For each  $j=1,\ldots,d$, let $m_j\geq 0$ be an integer and $\BFtheta_j = (\theta_{j,0},\theta_{j,1},\ldots,\theta_{j,m_j+1})\in\Real^{m_j+2}_{>0}$.
A GIBF covariance function is defined by
\begin{equation}
\begin{aligned}
K_{\mathsf{GIBF}}(\BFx,\BFx';\BFm,\BFtheta) ={}& \prod_{j=1}^d K_j(x_j,x_j'; m_j, \BFtheta_j), \quad \BFx,\BFx'\in\Real_{\geq 0}^d,\\[0.5ex]
K_j(x_j,  x'_j;m_j,\BFtheta_j) ={}& \sum_{\ell=0}^{m_j} \theta_{j, \ell}\frac{(x_j  x'_j)^\ell}{(\ell!)^2} + \theta_{j, m_j+1}
\int_0^\infty \frac{(x_j-u)_+^{m_j}( x'_j-u)_+^{m_j}}{(m_j!)^2}\, du,
\end{aligned}
\end{equation}
where $\BFm=(m_1,\ldots,m_d)$ and $\BFtheta = (\BFtheta_1,\ldots,\BFtheta_d)$, and  $(x)_+=\max(x,0)$ for $x\in\Real$.
\end{example}

A notable property, among others, of an $\BFm$-GIBF is that its sample path is $m_j$ times differentiable along the $j$-th coordinate, for each $j=1,\ldots,d$.
Moreover, if $m_j\geq 1$, then
the derivative of an $\BFm$-GIBF with respect to the $j$-th coordinate is a $(m_1,\ldots,m_j-1,\ldots,m_j+1,\ldots,m_d)$-GIBF.

\subsubsection{GP Regression}

If the simulation noise $\epsilon_{i,\ell}$ is independent of the GP prior and has a Gaussian distribution with a \emph{known} variance,
then the posterior distribution of $f$ is also a GP.
Let $\mathcal{D}_n\coloneqq \{(\BFx_i,\bar{y}_i):i=1,\ldots,n\}$  denote the simulation data.
The posterior mean function and covariance functions are given by
\begin{align}
\mu_n(\BFx) \coloneqq{}& \E[f(\BFx)|\mathcal{D}_n] = \mu(\BFx) + \BFk(\BFx)^\intercal(\BFK + \BFSigma)^{-1}(\bar{\BFy} - \BFmu), \label{eq:posterior-mean} \\
K_n(\BFx,\BFx')\coloneqq{}& \Cov[f(\BFx),f(\BFx')|\mathcal{D}_n] = K(\BFx,\BFx') - \BFk(\BFx)^\intercal (\BFK+\BFSigma)^{-1} \BFk(\BFx'), \label{eq:posterior-cov}
\end{align}
where $\BFk(\BFx)= (K(\BFx,\BFx_1),\ldots,K(\BFx,\BFx_n))^\intercal$, $\bar{\BFy}=(\bar{y}_1,\ldots,\bar{y}_n)^\intercal$, and
$\BFSigma$ is the $n$-by-$n$ diagonal matrix with the $i$-th diagonal element being $\sigma^2(\BFx_i)/r_i$.
Then, one can simply use $\mu_n(\BFx)$ to predict the response surface, that is,
\begin{equation}\label{eq:GPR-pred}
\hat f(\BFx) = \mu(\BFx) + \BFk(\BFx)^\intercal(\BFK + \BFSigma)^{-1}(\bar{\BFy} - \BFmu).
\end{equation}
Further, being the conditional expectation given the observations,
$\hat f(\BFx)$ is the best predictor that minimizes the mean squared prediction error; see \citet[page 153]{Rice07}.


\subsubsection{Selection of Hyperparameters}

A \emph{hyperparameter} is a parameter of a prior distribution in Bayesian statistics.
In the context of GP regression,
hyperparameters are those used to specify the mean function $\mu$ and the covariance function $K$.
For example, if we choose $\mu(\BFx) = \BFbeta^\intercal\BFphi(\BFx)$
and $K(\BFx,\BFx') = \tau^2 \exp(-\|\BFx-\BFx'\|^2/2\eta^2)$,
then the hyperparameters are $(\BFbeta, \tau, \eta)$.

Let $\BFtheta$ denote the collection of hyperparameters.
Let $\BFX$ denote the $n$-by-$d$ matrix whose $i$-th row is $\BFx_i^\intercal$ for all $i=1,\ldots,n$.
A usual approach to selecting $\BFtheta$ is to  maximize the following log ``likelihood'',
\begin{equation}\label{eq:log-likelihood}
\ln \mathsf{p}(\bar{\BFy}|\BFX) = -\frac{1}{2}(\bar{\BFy}-\BFmu(\BFtheta))^\intercal (\BFK(\BFtheta) + \BFSigma)^{-1}(\bar{\BFy}-\BFmu(\BFtheta)) - \frac{1}{2}\ln |\BFK(\BFtheta)| - \frac{n}{2}\ln (2\pi),
\end{equation}
where $|\cdot|$ denotes the determinant of a square matrix, and we write $ \BFmu(\BFtheta)$ and $\BFK(\BFtheta)$ to stress the dependence of $\BFmu$ and $\BFK$ on $\BFtheta$.
This is a nonlinear optimization problem with possibly multiple local optima.
Gradients of $\ln \mathsf{p}(\bar{\BFy}|\BFX)$ with respect to $\BFtheta$ can often be derived analytically and provided to
numerical optimization algorithms such as L-BFGS; see \cite{AnkenmanNelsonStaum10} for details.

A comment is warranted here regarding the notation of likelihood, however.
In Bayesian statistics, the likelihood---also termed the sampling distribution---is referring to the distribution of the observed data conditional on the data-generating process.
In the context of GP regression,
the likelihood is $\mathsf{p}(\bar{\BFy}|\BFX, \BFf)$,
where $\BFf=(f(\BFx_1),\ldots,f(\BFx_n))^\intercal$.
The conditioning on $\BFf$ is necessary
because $f$ is a \emph{random} function or sample path of the GP prior.
Note that $\ln \mathsf{p}(\bar{\BFy}|\BFX, \BFf)$ is \emph{not} identical to Equation~\eqref{eq:log-likelihood}.
Indeed,
\begin{align*}
\ln \mathsf{p}(\bar{\BFy}|\BFX, \BFf) ={}&  \ln \prod_{i=1}^n \frac{1}{\sigma(\BFx_i)\sqrt{2\pi}} \exp\biggl(-\frac{(\bar{y}_i - \mu(\BFx_i))^2}{2\sigma^2(\BFx_i)}\biggr) \\
={}&  -  \frac{1}{2}\sum_{i=1}^n (\bar{\BFy}-\BFmu(\BFtheta))^\intercal \BFSigma^{-1} (\bar{\BFy}-\BFmu(\BFtheta)) -\frac{n}{2}\ln(2\pi).
\end{align*}
Technically, $\mathsf{p}(\bar{\BFy}|\BFX)$ is called the log \emph{marginal likelihood} in Bayesian statistics, because it is the marginalization over $\BFf$:
\[
\mathsf{p}(\bar{\BFy}|\BFX) = \int \mathsf{p}(\bar{\BFy}|\BFX,\BFf) \mathsf{p}(\BFf|\BFX)\,d\BFf,
\]
where $\mathsf{p}(\BFf|\BFX)$ is the prior of $\BFf$,
which is a multivariate normal distribution; see \citet[Chapter 5]{RasmussenWilliams06} for more discussion.

\subsection{A Connection via Ridge Regularization}\label{sec:connect}

We start with the linear basis function model \eqref{eq:lienar-basis} for the case that many basis functions are included, even to the point of overparameterization $p \gg n$.
In this case, the OLS solution \eqref{eq:OLS-sol} is numerically unstable because the matrix $\BFPhi\BFPhi^\intercal $ is nearly singular if $p$ is smaller than but close to $n$ and becomes singular if $p\geq n$.
As a result,
$\|\hat\BFbeta\|$---the Euclidean norm of $\hat\BFbeta$---would explode, and
the prediction power of \eqref{eq:basis-pred} would be poor.

We now add ridge regularization \citep{Hastie20} to the least squares formulation to penalize the magnitude in norm of the solution, resulting in the method of regularized least squares (RLS):
\begin{equation}\label{eq:RLS}
\min_{\BFbeta} \frac{1}{n}\sum_{i=1}^n (\bar{y}_i - \BFbeta^\intercal \BFphi(\BFx_i))^2 + \lambda \|\BFbeta\|^2,
\end{equation}
where $\lambda \geq 0$ is the regularization parameter that controls the level of penalization.
The solution is
\begin{equation}\label{eq:RLS-solution}
\hat\BFbeta_\lambda = \BFPhi^\intercal (\BFPhi \BFPhi^\intercal + n\lambda \BFI )^{-1}\bar{\BFy},
\end{equation}
where $\BFI$ is the $n$-by-$n$ identity matrix.
The predictor for $f$ is then
\begin{equation}\label{eq:RLS-pred}
\hat f(\BFx) = \hat\BFbeta_\lambda^\intercal \BFphi(\BFx) = ( \BFPhi \BFphi(\BFx))^\intercal (\BFPhi \BFPhi^\intercal + n\lambda\BFI )^{-1}\bar{\BFy}.
\end{equation}

Note that
\[
\BFPhi \BFphi(\BFx)  = \begin{pmatrix}
\BFphi(\BFx)^\intercal \BFphi(\BFx_1) \\ \vdots \\ \BFphi(\BFx)^\intercal \BFphi(\BFx_n) \end{pmatrix}
\qq{and}
\BFPhi\BFPhi^\intercal = \begin{pmatrix}
\BFphi(\BFx_1)^\intercal \BFphi(\BFx_1)& \cdots & \BFphi(\BFx_1)^\intercal \BFphi(\BFx_n) \\
\vdots & \ddots & \vdots \\
\BFphi(\BFx_n)^\intercal \BFphi(\BFx_1)& \cdots & \BFphi(\BFx_n)^\intercal \BFphi(\BFx_n)\end{pmatrix}.
\]
Thus, the predictor \eqref{eq:RLS-pred} depends on the basis functions $\BFphi$ only through the product  of the form $\BFphi(\BFx)^\intercal\BFphi(\BFx')$ for some $\BFx$ and $\BFx'$.
If we define a bivariate function $K(\BFx,\BFx') = \BFphi(\BFx)^\intercal\BFphi(\BFx')$, then
the predictor \eqref{eq:RLS-pred} can be written as
\begin{equation}\label{eq:KRR-pred}
\hat f(\BFx) = \BFk(\BFx)^\intercal (\BFK + n\lambda\BFI)^{-1}\bar{\BFy},
\end{equation}
where $\BFk(\BFx)= (K(\BFx,\BFx_1),\ldots,K(\BFx,\BFx_n))^\intercal$
and $\BFK = (K(\BFx_i, \BFx_{i'}))_{i,i'=1}^n$.
Hence, the predictor \eqref{eq:RLS-pred} is formally identical to the GP regression predictor \eqref{eq:GPR-pred}, provided that $\mu\equiv 0$ and $\BFSigma = n\lambda \BFI$, i.e.,
$\sigma^2(\BFx_i)/r_i = n\lambda$ for all $i=1,\ldots,n$.

This connection implies that
the RLS method on linear basis function models may be interpreted, from a Bayesian perspective,
as GP regression.
That is, we impose  on $f$ a GP prior with mean $0$ and covariance function $K(\BFx,\BFx') = \BFphi(\BFx)^\intercal \BFphi(\BFx')$; moreover, the observation noise $\bar{\epsilon}_i$ is Gaussian with variance $n\lambda$.
In particular, with $\lambda=0$, we obtain an equivalence between the predictor \eqref{eq:basis-pred} and the noise-free GP regression, which is also known as GP interpolation.

The converse way of interpretation---GP regression as RLS---is also available, but it involves the theory of reproducing kernel Hilbert spaces, which is beyond the scope of this tutorial. We refer interested readers to \cite{KanHenSejSri18}.

\subsection{Enhancing Surrogates with Auxiliary Information}\label{sec:enhance}

So far we have treated the simulation model as a black box that performs nothing but transforming inputs to (noisy) outputs.
By doing so, we have implicitly assumed that the only data available when constructing a surrogate are simulation outputs.
However, there possibly exists auxiliary information in practice that we can leverage, without much extra computational overhead, to enhance the prediction capability of the surrogate.
We present below two general approaches.

\subsubsection{Enhancement with Stylized Models}
Simulation models, by design, are used to describe in detail the interactions between the components that constitute a complex stochastic system.
However, the main features of the same system may be captured by a stylized model that yields analytical expressions for the performance measure of interest as a function of the design variables,
provided that sufficiently many simplifying assumptions are made.

Consider the following example in \cite{ShenHongZhang18}.
The patient flow through various medical units of a hospital is a complicated queueing network.
The finite capacity of  one medical unit to accommodate patients often results in blocking patients from entering, making them stay in the upstream medical units even if the service has been completed there, and possibly creating further blocking.
A stylized model---however crude it may appear---is to decompose the network into isolated independent units by discarding the interactions and model each unit as a multi-server queue.
Performance measures such as mean length of stay can be derived in closed form for the stylized model.

Stylized models are often used to gain insights into the inner workings of the stochastic system of interest.
But they can be easily incorporated to construct a surrogate for the simulation model.
Let $\psi(\BFx)$ denote the analytical expression derived from a stylized model.
Then, we may simply add $\psi$ to the set of generic basis functions and use the augmented set either directly in linear basis function models
or to construct the mean function for GPs.

The key idea here is that one may use a crude---but computationally cheap otherwise---model to capture the overall trend of the response surface.
Conceivably,
the residual after de-trending will  have fewer variations, thereby easier to fit by a surrogate.
Thus, the requirement for $\psi$ to possess an analytical expression can be relaxed
as long as it can be computed sufficiently  fast,
e.g., from a low-fidelity simulation model. See \cite{LinMattaShanthikumar19} for more discussion.

\subsubsection{Enhancement with Gradient Observations}

Given simulation outputs for an input value $\BFx$,
an estimate of the gradient of the response surface with respect to $\BFx$ can often be obtained with a  negligible additional computational burden.
This kind of \emph{direct} gradient estimation---in contrast to finite-difference approximations---requires no re-simulation and
can be achieved via
infinitesimal perturbation analysis or the likelihood ratio/score function method  \citep{LEcuyer90}
in many simulation applications, including queueing systems \citep{Fu15} and financial engineering \citep[Chapter 7]{Glasserman03}.

Suppose that
in addition to $y_{i,\ell}$, the observation of $f(\BFx_i)$ for replication $\ell$ at point $\BFx_i$,
we obtain an unbiased estimate of  $\pdv{f(\BFx_i)}{x_j}$, the partial derivative of $f(\BFx_i)$ with respect to the $j$-th coordinate, and denote it by  $g_{i,j,\ell}$.

We first consider enhancing the linear basis function model with gradient observations.
We present below a formulation that generalizes the approach in \cite{FuQu14},  which focuses on linear regression models.
By doing so, it can be unified with the approach in \cite{ChenAnkenmanNelson13} that enhances GP surrogates with gradient observations.

Assume $f(\BFx) = \BFbeta^\intercal \BFphi(\BFx)$ and $\BFphi(\BFx)$ is differentiable.
Then, the gradient surface is $\nabla f(\BFx) = \BFbeta^\intercal \nabla \BFphi(\BFx)$.
Further, assume that for all $i=1,\ldots,n$ and $\ell=1,\ldots,r_i$,
\begin{equation}\label{eq:DiGAR}
\begin{aligned}
y_{i,\ell} ={}& \BFbeta^\intercal \BFphi(\BFx_i) + \epsilon_{i,\ell},  \\
g_{i,j,\ell} ={}& \BFbeta^\intercal
\pdv{\BFphi(\BFx_i)}{x_j} + \zeta_{i,j,\ell},  \quad j=1,\ldots,d,
\end{aligned}
\end{equation}
where the vector of noise terms $(\epsilon_{i,\ell}, \zeta_{i,1,\ell}, \ldots, \zeta_{i,d,\ell})$ has a multivariate normal distribution with mean vector $\BFzero$ and covariance matrix $\BFV$,
but no dependence exists across $i$ or $\ell$.
The correlation between, say, $\epsilon_{i,\ell}$ and
$\zeta_{i,j,\ell}$ stems from the fact that the gradient estimate $g_{i,j,\ell}$ is usually computed based on $y_{i,\ell}$.
Moreover, we do not consider the use of common random numbers here, which would introduce dependence across different design points.

Taking averages across replications in \eqref{eq:DiGAR} results in
\begin{equation}\label{eq:DiGAR-avg}
\begin{aligned}
\bar{y}_{i} ={}& \BFbeta^\intercal \BFphi(\BFx_i) + \bar{\epsilon}_{i},  \\
\bar{g}_{i,j} ={}& \BFbeta^\intercal
\pdv{\BFphi(\BFx_i)}{x_j} + \bar{\zeta}_{i,j},  \quad j=1,\ldots,d,
\end{aligned}
\end{equation}
where
$\bar{g}_{i,j} \coloneqq {r_i}^{-1}\sum_{\ell=1}^{r_i} g_{i,j,\ell}$
and $\bar{\zeta}_{i,j} \coloneqq {r_i}^{-1}\sum_{\ell=1}^{r_i} \zeta_{i,j,\ell} $.
Then, the system of Equations~\eqref{eq:DiGAR-avg} can be viewed as a linear basis function model,
with an augmented set of basis functions $\{\BFphi(\BFx), \pdv{\BFphi(\BFx)}{x_1},\ldots,\pdv{\BFphi(\BFx)}{x_d}\}$ and a vector of outputs $(\bar{y}_{i}, \bar{g}_{i,1},\ldots,\bar{g}_{i,d})$.
Due to the existence of correlations between the noise terms $\bar{\epsilon}_i, \bar{\zeta}_{i,1},\ldots,\bar{\zeta}_{i,d}$,
the OLS estimator of $\BFbeta$ is not statistically efficient; that is, there exists another estimator with a smaller variance.
Instead, it is recommended to use the method of generalized least squares (GLS):
\begin{align*}
\hat\BFbeta_{\mathsf{GLS}}={}&\argmin_{\BFbeta} (\bar{\BFy}_+ - \BFPhi_+  \BFbeta )^\intercal\BFV^{-1}(\bar{\BFy}_+ - \BFPhi_+ )\\
={}& (\BFPhi_+^\intercal \BFV^{-1} \BFPhi_+)^{-1} \BFPhi_+^\intercal \BFV^{-1} \bar{\BFy}_+,
\end{align*}
where
\[
\bar{\BFy}_+ =
\begin{pmatrix}
\bar{y}_1 \\
\bar{g}_{1,1} \\
\vdots \\
\bar{g}_{1,d} \\
\vdots \\
\bar{y}_n \\
\bar{g}_{n,1} \\
\vdots \\
\bar{g}_{n,d}
\end{pmatrix} \in\Real^{n(d+1)}
\qq{and}
\bar{\BFPhi}_+ =
\begin{pmatrix}
\BFphi(\BFx_1)^\intercal \\
\bigl(\pdv{\BFphi(\BFx_1)}{x_1}\bigr)^\intercal \\
\vdots \\
\bigl(\pdv{\BFphi(\BFx_1)}{x_d}\bigr)^\intercal \\
\vdots \\
\BFphi(\BFx_n)^\intercal \\
\bigl(\pdv{\BFphi(\BFx_n)}{x_1}\bigr)^\intercal \\
\vdots \\
\bigl(\pdv{\BFphi(\BFx_n)}{x_d}\bigr)^\intercal
\end{pmatrix} \in\Real^{n(d+1)\times p}.
\]

GP surrogates can also be enhanced with gradient observations. We briefly overview the approach proposed in \cite{ChenAnkenmanNelson13}, and refer to \cite{QuFu14} and \cite{HuoZhangZheng18} for further developments.

Assume $f(\BFx)=\BFbeta^\intercal \BFphi(\BFx) + \mathcal{M}(\BFx)$, where $\mathcal{M}(\BFx)$ is a zero-mean GP with  covariance function $K(\BFx,\BFx')$.
Then, the observations of the response surface and its gradient satisfy
\begin{equation}
\begin{aligned}
\bar{y}_i ={}& \BFbeta^\intercal \BFphi(\BFx_i) + \mathcal{M}(\BFx_i) + \bar{\epsilon}_i,   \\
\bar{g}_{i,j} ={}& \BFbeta^\intercal
\pdv{\BFphi(\BFx_i)}{x_j} + \pdv{\mathcal{M}(\BFx_i)}{x_j} + \bar{\zeta}_{i,j},\quad j=1,\ldots,d,
\end{aligned}
\end{equation}
where the partial derivative of a GP $\pdv{\mathcal{M}(\BFx)}{x_j}$  is defined in a mean-square sense.
Under regularity conditions,
it can be shown that
$(\mathcal{M}(\BFx), \pdv{\mathcal{M}(\BFx)}{x_1},\ldots,\pdv{\mathcal{M}(\BFx)}{x_d})$ forms a multi-output GP with mean zero,
and its covariance function can be derived explicitly by taking partial derivatives of $K(\BFx,\BFx')$.
The prediction can be made in closed form, but the formula is fairly involved so we omit the details.

\section{SO with Surrogates as Local Approximations}\label{sec:local}

Based on convergence guarantees, SO algorithms may be classified into three categories: locally convergent algorithms that converge to the set of local optimal solutions or stationary points, globally convergent algorithms that converge to the set of global optimal solutions, and heuristic algorithms that have no convergence guarantee. Although global convergence is ideal, it is a global property, i.e., it typically requires exploring the entire feasible region in the limit. We will introduce many such algorithms in Section~\ref{sec:global}.

In many practical situations, however, the computational budget is limited and only allows exploring a small proportion of the feasible region. In these situations, local convergence that only requires local information becomes more meaningful, and it may be practically tested with a certain statistical guarantee; see, for instance, \cite{BettonvilCastilloKleijnen09} and \cite{XuNelsonHong10}.
One way to obtain the local information about a solution is through a local surrogate, which allows the SO algorithm to check whether the solution is a local optimal solution and, if not, identify a direction (or a region) where better solutions may be found.

For SO purposes, due to Taylor’s expansion, the most natural choices of local surrogates are low-order polynomials (see Section~\ref{sec:polynomial}), especially first- and second-order polynomials. Response-surface methodology (RSM) is a collection of statistical methods that build on this idea to solve stochastic optimization problems, which include SO problems. We review RSM in Section~\ref{sec:RSM}. However, RSM often requires human involvement because it typically deals with real experiments (such as agricultural or clinical experiments), so it is not particularly suited for SO problems. The stochastic trust-region response-surface method (STRONG) of \cite{ChangHongWan13} solves this problem by combing RSM with the trust-region method developed for deterministic nonlinear optimization. We review the STRONG and the related algorithms in Section~\ref{sec:STRONG}. Other than low-order polynomials, other types of local surrogates have also been used in SO. We introduce the surrogate-based promising area search algorithm of \cite{FanHu18} in Section~\ref{sec:SPAS}.

\subsection{Response Surface Methodology}\label{sec:RSM}

RSM was first developed by \cite{BoxWilson51} to optimize the operating conditions of a chemical process. It has evolved into a major tool for optimizing real (i.e., non-simulation) experiments. According to
\cite{RSM_textbook}, ``RSM is a collection of statistical and mathematical techniques useful for developing, improving, and optimizing processes.” It typically includes two stages. In the first stage, it runs a number of experiments in a local region of the current solution and builds a first-order surrogate in the form of
\begin{equation}\label{eq:3_linear}
f({\BFx}) = \beta_0+\sum_{j=1}^d \beta_j x_j.
\end{equation}
Notice that Equation (\ref{eq:3_linear}) implies that $\nabla f(\BFx)=(\beta_1,\ldots,\beta_d)^\intercal$. The surrogate essentially provides an ascent direction to allow RSM to find a better solution. Once it has a new solution, it repeats the process to find a better solution iteratively until the first-order model is no longer adequate. Then, RSM switches to the second stage, where it builds a second-order surrogate in the form of Equation~\eqref{eq:quadratic} to locate the optimal solution.

Because RSM is typically used for expensive real experiments (or simulation experiments that are slow to run), large-sample properties, e.g., convergence or rate of convergence, are typically not considered in the literature, and the focus is mainly on statistical issues, such as the design of experiments (DOE) to efficiently estimate the surrogates and the tests of model inadequacy and optimality \citep{RSM_textbook}. For instance, there are $d+1$ and $d(d+1)/2+1$ parameters in the first- and second-order models, respectively. Different DOE schemes are proposed to reduce the required number of experiments to appropriately estimate these parameters \citep{Kleijnen15}.

RSM has also become a popular heuristic tool for SO  \citep{HoodWelch93,Kleijnen08}, especially when the simulation experiments are time-consuming. Much has been developed to understand how simulation experiments impact the statistical properties of RSM. For instance, \cite{SchrubenMargolin78}  study how the use of common random numbers impacts the fitting of polynomials;  \cite{AngunKleijnenHertogGurkan09} consider stochastic constraints; and \cite{BettonvilCastilloKleijnen09} develop tests of the Karush–Kuhn–Tucker conditions.

\subsection{Stochastic Trust-Region Response-Surface Method}\label{sec:STRONG}

While RSM is a popular tool for SO, it has several problems, especially when simulation experiments are relatively fast (so that a large number of experiments may be conducted) and simulation noise is significant. First, it requires human involvement. For instance, in each iteration of the RSM, a surrogate needs to be optimized in a local region, but the local regions are determined by experimenters based on their experience. Moreover, the transitions between first- and second-order models typically also depend on human experience, and it is not clear whether a second-order model may transition back to a first-order model if it is found inadequate due to the simulation noise. Second, it is not clear whether the RSM algorithms have any convergence guarantees. This is a relevant question, especially when the number of simulation experiments becomes large.

\cite{ChangHongWan13}  propose the STRONG algorithm, which combines the trust-region method of deterministic nonlinear optimization with the RSM framework, to solve the two problems of RSM. In any iteration, say iteration $k$, let ${\BFx}_k$ and $\Delta_k$ denote the current solution and the size of the trust region. STRONG conducts the following four steps:
\begin{enumerate}[label=\emph{Step \arabic*.}, left=0pt]
\item Construct a local model $r_k({\BFx})$ around the current solution ${\BFx}$. If $\Delta_k\ge\tilde\Delta$, where $\tilde\Delta$ is a threshold, $r_k({\BFx})$ is a first-order model; otherwise, $r_k({\BFx})$ is a second-order model;
\item Solve ${\BFx}_k^*\in{\argmax}\{r_k({\BFx}): {\BFx}\in{\cal B}({\BFx}_k,\Delta_k)\}$, where ${\cal B}({\BFx}_k,\Delta_k)$ denotes a $d$-dimensional ball centered at ${\BFx}_k$ with a radius $\Delta_k$, which is the trust region at iteration $k$;
\item Simulate a number of observations at ${\BFx}_k^*$ and estimate $f({\BFx}_k^*)$;
\item Conduct the sufficient-reduction and ratio-comparison tests to update ${\BFx}_{k+1}$ and $\Delta_k$.
\end{enumerate}

In the algorithm, the sufficient-reduction test conducts new simulation experiments to test whether $\BFx_k^*$ is statistically significantly better than ${\BFx}_k$. If it is not, then the current solution is not updated, i.e., ${\BFx}_{k+1}={\BFx}_k$, and the trust region shrinks, i.e., $\Delta_{k+1}=\gamma_1\Delta_k$ where $0<\gamma_1<1$ is a constant. If $\BFx_k^*$ passes the sufficient-reduction test, then the algorithm moves to the ratio-comparison test, which computes
\[\rho_k = \frac{\bar f_k({\BFx}_k^*) -\bar f_k({\BFx}_k)}{ r_k({\BFx}_k^*) - r_k({\BFx}_k)},
\]
where $\bar f_k({\BFx}_k^*)$ and $\bar f_k({\BFx}_k)$ are the estimated objective values (using simulation experiments) at ${\BFx}_k^*$ and ${\BFx}_k$, respectively. Notice that $\rho_k$ denotes the ratio between the actual observed improvement and the predicted improvement. Let $0<\eta_0<\eta_1<1$ be two thresholds. If $\rho_k\ge\eta_1$, which implies that the local model works well, the algorithm then moves the current solution to the new solution, i.e., ${\BFx}_{k+1}={\BFx}_k^*$, and enlarges the size of the trust region, i.e., $\Delta_{k+1}=\gamma_2\Delta_k$ where $\gamma_2>1$ is a constant. If $\eta_0\le\rho_k< \eta_1$, which implies that the local model has some predictive power, the algorithm updates the new solution i.e., ${\BFx}_{k+1}={\BFx}_k^*$, but keeps the size of the trust region, i.e., $\Delta_{k+1}=\Delta_k$. If $\rho_k<\eta_1$, which implies that the local model works poorly, the algorithm keeps the current solution, i.e., ${\BFx}_{k+1}={\BFx}_k$, and shrinks the size of the trust region, i.e., $\Delta_{k+1}=\gamma_1\Delta_k$.

Notice that the STRONG algorithm uses the trust region as the local region and its size $\Delta_k$ is adaptively updated based on the sufficient-reduction and ratio-comparison tests. Furthermore, the transitions between first- and second-order models are based on the size of the trust region $\Delta_k$. If it is larger than the threshold $\tilde\Delta$, a first-order model is used; otherwise, a second-order model is used. This transition rule allows the algorithm to transition between the two models in both directions. Therefore, the algorithm gets rid of the human involvement that is necessary for typical RSM algorithms. Furthermore,  \cite{ChangHongWan13} show that, under certainly technical conditions on the estimated surrogates, the STRONG algorithm converges to a stationary point of the original SO problem.

The use of the trust-region method in SO has also been studied by others. For instance, \cite{DengFerris09} combine it with the sample-average approximation to solve SO problems; \cite{ShashaaniHashemiPasupathy18} integrate it into a derivative-free algorithm to solve SO problems; and \cite{MathesenPedrielliNg17} use it with a restart approach to design a globally convergent SO algorithm. The STRONG algorithm has also been extended by \cite{ChangLiWan14} to include a screening stage so that it may solve large-scale SO problems with hundreds of dimensions.

\subsection{Surrogate-Based Promising Area Search}\label{sec:SPAS}

In Sections~\ref{sec:RSM} and \ref{sec:STRONG} we introduced the RSM and STRONG algorithms. Both of them build low-order polynomials as local surrogates and use these models to guide the optimization process. In this subsection we introduce the Surrogate-Based Promising Area Search (SPAS) algorithm of \cite{FanHu18} that allows the use of any interpolation surrogates, e.g., kriging, splines or radial basis functions.

Unlike RSM and STRONG, SPAS fits a global surrogate in each iteration. However, as \cite{FanHu18} pointed out, ``the use of the surrogate in our approach [i.e., SPAS] is not intended to provide a global fit of the underlying response surface, but rather aims to accurately predict the objective function values at unsampled points within the current search area.'' That’s why we also include SPAS in Section~\ref{sec:local}, which focuses on SO algorithms with local surrogate approximations.

In each iteration, SPAS consists of the following four steps:
\begin{enumerate}[label=\emph{Step \arabic*.}, left=0pt]
\item Construct the most promising area (MPA), and sample a set of candidate solutions from it;
\item Estimate the objective values of all visited solutions using a shrinking-ball method, which averages the samples in a $d$-dimensional ball centered at the solution;
\item Build a surrogate that interpolates all these solutions;
\item Find the optimal solution of the surrogate within the MPA.
\end{enumerate}

Besides the use of surrogates, SPAS  also integrates several other critical ideas of SO. The shrinking-ball method of estimating the function value at any solution was first introduced by \cite{BaumertSmith02}, and it has been studied and applied by \cite{AndradottirPrudius10} and \cite{KiatsupaibulSmithZabinsky18}. The concept of the MPA was first proposed by \cite{HongNelson06} in their COMPASS algorithm, which solves DOvS problems. \cite{HongNelson07} further extend the idea into a general DOvS framework. By combining surrogate modeling, the shrinking-ball method and the MPA, SPAS is proved to converge to the set of local optimal solutions if the objective function is Lipschitz continuous.

\section{SO with Surrogates as Global Approximations}\label{sec:global}

There are mainly two strategies for using global surrogates such as linear basis function models and GPs to solve a continuous SO problem, depending on whether the design points  $\mathcal{X} = \{\BFx_1,\ldots,\BFx_n\}$ are chosen in a \emph{static} fashion or a \emph{sequential} one.

The former means that $\mathcal{X}$ is determined---once and for all---prior to any simulation experiments.
With no observations of the response surface being available, a primary design principle of $\mathcal{X}$ is to cover the design space as much as possible
so that most part of the response surface can be recovered after the observations are obtained.
Typical experimental designs include lattice designs and space-filling designs,
and we refer to \citet[Chapter 5]{SantnerWilliamsNotz03} for details.
Given $\mathcal{X}$,
one runs simulation at each design point, possibly multiple times, fits a surrogate with the observations, and then optimizes the predicted surface $\hat f(\BFx)$ induced by the surrogate.
Being a deterministic function, $\hat f(\BFx)$ can be optimized with
any numerical optimization algorithms \citep{NocedalWright06}.
We refer to \cite{BartonMeckesheimer06} for more discussion on
the use of global surrogates in conjunction with static experimental designs.

Recent advances in SO with surrogates as global approximations are dominated by sequential experimental designs---design points are determined one at a time after each new observation of the response surface is made.
Each new design point is selected based on
(i) the updated surrogate reflecting the previous observations, and
(ii) certain criteria that balance \emph{exploration} and \emph{exploitation}.
By exploration, we mean searching the part of the design space that has much uncertainty;
by exploitation, we mean searching the area in the proximity of the current best solution.
The need for quantifying the uncertainty about the response surface renders GPs the most popular class of surrogates.
(Recall that posterior distributions of a GP can be derived in closed form and the computation is reduced to linear algebra.)
This line of research is closely related to Bayesian optimization \citep{ShahriariSwerskyWangAdamsdeFreitas16,Frazier18}
in the machine learning literature,
where a primary motivation is hyperparameter tuning of sophisticated machine learning algorithms/models \citep{FeurerHutter19}.
We introduce below three representative examples of such methods:
knowledge gradient \citep{ScottFrazierPowell11},
upper confidence bound \citep{SrinivasKrauseKakadeSeeger12},
and GP-based search \citep{SunHongHu14}.
All three methods share the following structure procedure-wise:
\begin{enumerate}[label=\emph{Step \arabic*.}, ref=\arabic*, left=0pt]
\item Impose a GP prior on $f$;
\item Select the next batch of design points subject to a prescribed ``criterion'' that is computed using the current belief about $f$;\label{step:selec}
\item Run simulation experiments at each of the newly selected design points;
\item Update the GP posterior given the new observations of $f$ via Equations~\eqref{eq:posterior-mean} and \eqref{eq:posterior-cov};\label{step:update}
\item Repeat Steps~\ref{step:selec}--\ref{step:update} until the simulation budget is exhausted;
\item Optimize the posterior mean function and return the optimum as a solution to problem~\eqref{eq:SO}.
\end{enumerate}

As demonstrated below, GP-based sequential methods for continuous SO problems mainly differ in how to define the criterion in Step~\ref{step:selec} for selecting new design points.


\subsection{Knowledge Gradient}
Knowledge gradient (KG) was originally proposed to solve R\&S problems \citep{FrazierPowellDayanik08,FrazierPowellDayanik09}.
The method was generalized in \cite{ScottFrazierPowell11} to cover continuous SO problems.

Suppose that simulation experiments have been made at $\{\BFx_1,\ldots,\BFx_n\}$ with one replication each,
generating observations $y_i = f(\BFx_i) + \epsilon_i$,
$i=1,\ldots,n$.
We are interested in selecting the next design point $\BFx_{n+1}$.
Let $\mathcal{D}_n= \{(\BFx_i,y_i):i=1,\ldots,n\}$.
Let $\mu_n$ and $K_n$ be the posterior mean and covariance functions of $f$ conditional on $\mathcal{D}$.
Following Equations~\eqref{eq:posterior-mean} and \eqref{eq:posterior-cov},
it can be shown that $\mu_n$ and $K_n$ satisfy the following updating scheme:
\begin{align}
\mu_{n+1}(\BFx) ={}&  \mu_n(\BFx) +  \delta_n(\BFx, \BFx_{n+1}) Z_{n+1}, \label{eq:updating-mean} \\
K_{n+1}(\BFx,\BFx') = {}&  K_n(\BFx,\BFx') - \delta_n(\BFx, \BFx_{n+1})\delta_n(\BFx', \BFx_{n+1}), \label{eq:updating-cov}
\end{align}
where
\[
Z_{n+1} \coloneqq \frac{y_{n+1}-\mu_n(\BFx_{n+1})}{\sqrt{K_n(\BFx_{n+1}, \BFx_{n+1}) + \sigma^2(\BFx_{n+1})}}
\qq{and}
\delta_n(\BFx, \BFv)\coloneqq \frac{K_n(\BFx,\BFv)}{\sqrt{K_n(\BFv,\BFv)+ \sigma^2(\BFv)}};
\]
moreover, $Z_{n+1}$ is a standard normal random variable conditional on $\mathcal{D}_n$.

The KG method selects $\BFx_{n+1} = \argmax_{\BFx\in\mathscr{X}} \mathrm{KG}_n(\BFx)$, where
\begin{equation}\label{eq:KG}
\mathrm{KG}_n(\BFx)\coloneqq \E\Bigl[\max_{\BFu\in\mathscr{X}}\mu_{n+1}(\BFu) -  \max_{\BFu\in\mathscr{X}}\mu_n(\BFu) \big|\mathcal{D}_n, \BFx_{n+1}=\BFx\Bigr].
\end{equation}
The interpretation of $\mathrm{KG}_n(\BFx)$ is as follows.
If our simulation budget were exhausted after collecting data $\mathcal{D}_n$, then we would
use $\max_{\BFu}\mu_n(\BFu)$ to estimate the maximum value of $f$.
However, now that we are allowed to run one more simulation experiment at $\BFx_{n+1}$,
the posterior mean function will become $\mu_{n+1}(\cdot)$ in Equation~\eqref{eq:updating-mean} after the new data point $(\BFx_{n+1},y_{n+1})$ is obtained.
Thus, the increment in the estimated maximum value of $f$ is $\max_{\BFu}\mu_{n+1}(\BFu) -  \max_{\BFu}\mu_n(\BFu)$.
Before the simulation is run at $\BFx_{n+1}$, this increment is a random variable conditional on $\mathcal{D}_n$, and its distribution is determined by the standard normal random variable $Z_{n+1}$ in Equation~\eqref{eq:updating-mean}.

There are several approaches for the numerical maximization of $\mathrm{KG}_n(\BFx)$.
\cite{ScottFrazierPowell11} propose  maximizing
\[
\widetilde{\mathrm{KG}}_n(\BFx)\coloneqq \E\Bigl[\max_{1\leq i\leq n+1}\mu_{n+1}(\BFx_i) -  \max_{1\leq i \leq n+1}\mu_n(\BFx_i) \big|\mathcal{D}_n, \BFx_{n+1}=\BFx \Bigr],
\]
a discrete proxy of $\mathrm{KG}_n(\BFx)$,
because $\widetilde{\mathrm{KG}}_n(\BFx)$ and
its gradient with respect to $\BFx$
can both be computed explicitly.

A second approach to solving  $\max_{\BFx}\mathrm{KG}_n(\BFx)$ is to view it as a stochastic optimization problem,
in which the only random variable involved is the standard normal $Z_{n+1}$.
Then, one may apply sample average approximation \citep{KimPasupathyHenderson15} by simulating realizations of $Z_{n+1}$.

Yet another approach, proposed in \cite{WuFrazier16}, is to apply stochastic approximation \citep{ChauFu15}.
Note that, under mild regularity conditions,
\[
\nabla_{\BFx} \mathrm{KG}_n(\BFx) = \nabla_{\BFx}
\E\Bigl[\max_{\BFu\in\mathscr{X}}\mu_{n+1}(\BFu) -  \max_{\BFu\in\mathscr{X}}\mu_n(\BFu) \big|\mathcal{D}_n, \BFx_{n+1}=\BFx\Bigr]
= \E\Bigl[\nabla_{\BFx}\max_{\BFu\in\mathscr{X}}\mu_{n+1}(\BFu) \big|\mathcal{D}_n, \BFx_{n+1}=\BFx\Bigr].
\]
Hence, $\nabla_{\BFx}\max_{\BFu\in\mathscr{X}}\mu_{n+1}(\BFu)$ is an unbiased estimator of $\nabla_{\BFx} \mathrm{KG}_n(\BFx)$,
and it can be computed by applying the envelope theorem \citep{MilgromSegal02}.

\subsection{Upper Confidence Bound}
Upper confidence bound (UCB) is a celebrated class of methods for multi-armed bandit (MAB) problems.
Similar to R\&S problems, MAB problems are also concerned with finding the optimal among a finite set of alternatives with unknown performances/rewards.
A key difference between the two classes of problems lies in the objective.
Basically,
MAB can be viewed as an online decision-making problem
that aims to maximize the cumulative rewards collected over the entire time horizon;
in contrast,
R\&S is more of an offline flavor and focuses on the final reward collected at the end of the time horizon.
We refer to  \cite{Slivkins19} and \cite{Auer02} for introductions to MAB problems and UCB-type algorithms, respectively.

\cite{SrinivasKrauseKakadeSeeger12} generalize UCB to the setting of optimizing a GP sample path.
The general structure of the GP-UCB method is basically the same as that of the KG method,
except that the next design point is selected as  $\BFx_{n+1} = \argmax_{\BFx\in\mathscr{X}} \mathrm{UCB}_n(\BFx)$, where
\begin{equation}\label{eq:GPUCB}
\mathrm{UCB}_n(\BFx)\coloneqq \mu_n(\BFx) +  \sqrt{\gamma_n K_n(\BFx,\BFx)},
\end{equation}
and $\gamma_n>0$ is a tuning parameter that should grow as a function of $n$.
The form of $\mathrm{UCB}_n(\BFx)$ clearly shows  the trade-off between exploration and exploitation.
A potential design point $\BFx$ is favored if either $\mu_n(\BFx)$ is large (exploitation), or $K_n(\BFx,\BFx)$ is large (exploration).

Evidently, $\mathrm{UCB}_n(\BFx)$ is much easier to maximize than $\mathrm{KG}_n(\BFx)$,
as the former involves no expectation.
Nevertheless, a potential downside of the GP-UCB method is that its performance depends critically on the choice of $\gamma_n$.
It might be tempting to make $\gamma_n$ grow at a rate of $\ln(n)$,
which is both a typical choice for MAB problems \citep{Auer02}
and the choice analyzed in \cite{SrinivasKrauseKakadeSeeger12}.
However, such a choice is recommended with the aim of maximizing the cumulative reward of the form $\sum_{i=1}^n f(\BFx_i)$.
There may conceivably exist a better guideline for setting  $\gamma_n$ for problem~\eqref{eq:SO}.
This issue is yet to be addressed.

\subsection{GP-Based Search}

Random search is a general class of algorithms for SO, which samples randomly a number of design points in each iteration of the algorithm based a sampling distribution that may be updated based on all information collected through the optimization process. Random search algorithms are most popular for DOvS problems, but they have also been developed to solve continuous SO problems; see, for instance, the recent reviews of \cite{HongNelsonXu15} and \cite{Andradottir15}. One of the critical issues in developing random search algorithms is how to design a sampling distribution that automatically balances exploration and exploitation. As pointed out earlier in this section, GP is capable of quantifying the uncertainty of the response surface. Therefore, it can be used to facilitate the design of good sampling distributions, which is the basic idea of the GP-based Search (GPS) algorithm of \cite{SunHongHu14}.

\cite{SunHongHu14} note that the posterior mean function $\mu_n(\BFx)$ and the posterior variance function $\sigma_n^2(\BFx)=K_n(\BFx,\BFx)$ of Equations (\ref{eq:posterior-mean}) and (\ref{eq:posterior-cov}) provide information on exploitation and exploration, respectively. In particular, higher mean values indicate the region needs more exploitation and higher variance values indicate the region needs more exploration. Therefore, the GPS algorithm uses the following sampling distribution
\[
h(\BFx) = \frac{\Pr{Z(\BFx)>c}}{\sum_{{\BFz}\in\mathscr{X}} \Pr{Z({\BFz})>c}},
\]
where $\mathscr{X}$ is a finite set of discrete solutions, $c$ is set as the current estimated optimal value, and $Z(\BFx)$ follows a normal distribution with mean $\mu_n(\BFx)$ and variance $\sigma_n^2(\BFx)$ for any $\BFx\in\mathscr{X}$. Notice that the sampling distribution combines both the mean and variance information and balances exploration and exploitation seamlessly.

To use the sampling distribution there are two remaining issues. The first is how to sample from the distribution. Notice that the denominator of $h(\BFx)$ involves a summation that is typically difficult to compute. \cite{SunHongHu14} solve the problem by developing an acceptance-rejection algorithm and a Markov chain Monte Carlo algorithm to sample from the distribution.

The second issue is the calculation of $\mu_n(\BFx)$ and $\sigma_n^2(\BFx)$ using Equations (\ref{eq:posterior-mean}) and (\ref{eq:posterior-cov}). Notice that the calculation involves a matrix inversion. When the number of design points is large, this calculation is time-consuming. Furthermore, when the design points are close to each other (which is common in later iterations of the algorithm when it focuses more on exploiting the good regions), the matrix is often ill-conditioned and the inversion becomes difficult. To solve the problem, the GPS algorithm takes a very pragmatic view towards the GP. Instead of considering the objective value as a sample path from the GP, as in Bayesian optimization algorithms, the GPS algorithm only treats it as a surrogate approximation that facilitates the generation of good sampling distributions. It proposes the following GP to model the objective function:
\begin{equation}\label{eqn:fastkriging}
f({\BFx})=\mathcal{M}(\BFx)+\lambda(\BFx)^\intercal(\bar\BFy-{\BFM})+\lambda(\BFx)^\intercal{\cal E},
\end{equation}
where $\mathcal{M}(\BFx)$ is an unconditional GP, $\lambda(\BFx)=(\lambda_1(\BFx),\ldots,\lambda_n(\BFx))^\intercal$ is a vector of weight functions, ${\BFM}=(M(\BFx_1),\ldots,M(\BFx_n))^\intercal$ is a vector of $\mathcal{M}(\BFx)$ evaluated at $\BFx_1,\ldots,\BFx_n$ and ${\cal E}=(\epsilon_1,\ldots,\epsilon_n)^\intercal$ is an $n$-dimensional random vector following a multivariate normal distribution with the mean $\BFzero$ and covariance matrix $\BFSigma$, which is the $n$-by-$n$ diagonal matrix with the $i$-th diagonal element being $\sigma^2(\BFx_i)/r_i$. Furthermore, in Equation (\ref{eqn:fastkriging}), $\mathcal{M}(\BFx)$ and $\cal E$ are independent of each other and $\bar\BFy$ is considered deterministic when building the model.

Let $\tilde\mu_n(\BFx)$ and $\tilde\sigma_n^2(\BFx)$ denote the mean and variance functions of the new GP model. When the weight function vector $\lambda(\BFx)$ is continuous in $\BFx$ and it satisfies $\lambda_i(\BFx)\ge 0$, $\sum_{i=1}^n \lambda_i(\BFx)=1$ and $\lambda_i(\BFx_j)=1\{\BFx_i=\BFx_j\}$, where $1\{\cdot\}$ is the indicator function, \cite{SunHongHu14} show that
\begin{eqnarray*}
\tilde\mu_n(\BFx) &=& \lambda(\BFx)^\intercal\bar\BFy,\\
\tilde\sigma_n^2(\BFx) &=& K(\BFx,\BFx)-2\lambda(\BFx)^\intercal\BFk(\BFx)+\lambda(\BFx)^\intercal(\BFK+\BFSigma)\lambda(\BFx).
\end{eqnarray*}
Then, both $\tilde\mu_n(\BFx)$ and $\tilde\sigma_n^2(\BFx)$ may be calculated directly without matrix inversion. Furthermore, \cite{SunHongHu14} show that $\tilde\mu_n(\BFx)$ interpolates all design points through their sample means, i.e., $\tilde\mu_n(\BFx_i)=\bar y_i$ for all $i=1,\ldots,n$, and $\tilde\mu_n(\BFx)$ and $\tilde\sigma_n^2(\BFx)$ capture the information of exploitation and exploration and, therefore, can be used to construct the sampling distribution in the GPS algorithm.

The GPS algorithm has global convergence \citep{SunHongHu14}. However, it is designed to solve DOvS problems. \cite{SunHuHong18} extend it to solve continuous SO problems and proved its global convergence. To further simplify the sampling process, \cite{SunHuHong18} propose to approximate the sampling distribution by a Gaussian mixture model that can be sampled easily, and show that the resulting algorithm is still globally convergent.

\section{Computation for Large Datasets}\label{sec:largescale}

We have taken it for granted thus far in this tutorial that surrogates are computationally fast,
and rightly so in light of their analytical tractability.
However, this postulate is challenged when the number of design points is large
because it usually involves numerically inverting a large matrix to process the simulation data $\{(\BFx_i, \bar{y}_i):i=1,\ldots,n\}$; see, e.g., Equation~\eqref{eq:basis-pred} for linear basis function models, as well as Equations~\eqref{eq:posterior-mean} and \eqref{eq:posterior-cov} for GPs.
It is well known that the time complexity for matrix inversion scales as $\mathcal{O}(n^3)$ in general.
As $n$ grows, surrogates are increasingly demanding in computation and eventually become even more expensive than the simulation model that they aim to approximate in the first place \citep{HuangAllenNotzZeng06,SunHongHu14,SalemiSongNelsonStaum19}.

The need for processing large datasets calls for approximation methods to reduce the computational burden caused by numerical inversion of large matrices.
There exists a huge literature on approximate computation for GP regression; see \cite{LiuOngShenCai20} for a recent survey.
We present two popular methods---the Nystr\"om method and random features.
Both methods have drawn substantial interest in recent years.

\subsection{The Nystr\"om Method}
A central idea to mitigate the challenge of computing  $(\BFK+\BFSigma)^{-1}$ is to find a low-rank approximation of $\BFK$.
In particular, consider a rank-$m$ matrix of the form $\tilde{\BFK}=\BFU \BFC\BFV$, where $\BFU\in\Real^{n\times m}$, $\BFC \in\Real^{m\times m}$, and $\BFV\in\Real^{m\times n}$ with $m<n$.
If $\BFC$ is invertible,
the Woodbury matrix identity \citep[page 19]{HornJohnson12} asserts that
\begin{equation}\label{eq:woodbury}
(\BFK+\BFSigma)^{-1}\approx (\tilde{\BFK}+\BFSigma)^{-1} = \BFSigma^{-1} - \BFSigma^{-1}\BFU(\BFC^{-1} + \BFV\BFSigma^{-1}\BFU)^{-1}\BFV\BFSigma^{-1}.
\end{equation}
Then, the computational bottleneck has been transformed to the inversion of smaller, $m$-by-$m$ matrices. (The inversion of $\BFSigma$ is easy because it is a diagonal matrix.)

For ease of presentation, let $I=\{1,\ldots,n\}$ and $A\subset I$ be a subset of size $m$, which is also called the active set of indices.
The Nystr\"om method was originally devised to approximate the eigenfunctions of a covariance function $K(\cdot,\cdot)$; see, e.g., \citet{WilliamsSeeger01}.
It suggests the following low-rank approximation of the covariance matrix $\BFK$:
\begin{equation}\label{eq:nystrom}
\tilde{\BFK} = \BFK_{n,m} \BFK_{m,m}^{-1} \BFK_{m,n},
\end{equation}
where $\BFK_{n,m} \coloneqq (K(\BFx_i, \BFx_{i'}))_{i\in I, i'\in A}$,
$\BFK_{m,m} \coloneqq (K(\BFx_i, \BFx_{i'}))_{i\in A, i'\in A}$,
and $\BFK_{m,n} \coloneqq (K(\BFx_i, \BFx_{i'}))_{i\in A, i'\in I}$.
Then, the posterior mean function in  Equation~\eqref{eq:posterior-mean} can be approximated by replacing $\BFK$ with $\tilde{\BFK}$:
\begin{align}
&\E[f(\BFx)|\mathcal{D}_n] \nonumber \\
\approx{}& \mu(\BFx) + \BFk(\BFx)^\intercal(\tilde{\BFK}+ \BFSigma)^{-1}(\bar{\BFy} - \BFmu) \nonumber\\
={}& \mu(\BFx) + \BFk(\BFx)^\intercal \BFSigma^{-1} (\bar{\BFy} - \BFmu) -
\BFk(\BFx)^\intercal \BFSigma^{-1}\BFK_{n,m}(\BFK_{m,n} + \BFK_{m,n}\BFSigma^{-1}\BFK_{n,m})^{-1}\BFK_{m,n}\BFSigma^{-1} (\bar{\BFy} - \BFmu), \label{eq:posterior-mean-approx}
\end{align}
where the last step follows from \eqref{eq:woodbury}.
Despite its long expression, the approximation \eqref{eq:posterior-mean-approx} can be computed in time $\mathcal{O}(m^2n)$.
In addition, the posterior covariance function in Equation~\eqref{eq:posterior-cov} can be approximated in a similar fashion: \begin{align}
& \Cov[f(\BFx), f(\BFx')|\mathcal{D}_n] \nonumber  \\
\approx {}& K(\BFx,\BFx') - \BFk(\BFx)^\intercal (\tilde{\BFK}+\BFSigma)^{-1} \BFk(\BFx') \nonumber \\
={}&  K(\BFx,\BFx') - \BFk(\BFx)^\intercal \BFSigma^{-1} \BFk(\BFx') +
\BFk(\BFx)^\intercal \BFSigma^{-1}\BFK_{n,m}(\BFK_{m,n} + \BFK_{m,n}\BFSigma^{-1}\BFK_{n,m})^{-1}\BFK_{m,n}\BFSigma^{-1} \BFk(\BFx'), \label{eq:posterior-cov-approx}
\end{align}
whose time complexity is also $\mathcal{O}(m^2n)$.

However, there is a caveat in the above use of the low-rank approximation \eqref{eq:nystrom}.
That is, using \eqref{eq:posterior-cov-approx} to approximate the posterior variance $\Var[f(\BFx)|\mathcal{D}_n] = \Cov[f(\BFx), f(\BFx)|\mathcal{D}_n]$ may yield a negative value!
A better implementation of the Nystr\"om method is to construct a covariance function $\tilde{K}(\cdot,\cdot)$ to systematically replace the occurrences of $K(\cdot,\cdot)$.

Specifically, define $\tilde{K}(\BFx,\BFx') \coloneqq  \BFk_m(\BFx)^\intercal \BFK_{m,m}^{-1} \BFk_m(\BFx')$,
where $\BFk_m(\BFx) \in \Real^m$ is the vector composed of $K(\BFx,\BFx_i)$ for all $i\in A$.
It is easy to show that
(i) $\tilde{K}$ is a covariance function, and (ii)
the covariance matrix associated with evaluating $\tilde{K}$ at $\{\BFx_1,\ldots,\BFx_n\}$ is identical to $\tilde{\BFK}$ in Equation~\eqref{eq:nystrom}.
Let $\tilde{f}$ denote a GP with mean function $\mu$ and covariance function $\tilde{K}$.
Then, we can use the posterior distribution of $\tilde{f}$ to  approximate that of $f$.
This treatment is adopted by \cite{SmolaScholkopf00} and \cite{RudiCamorianoRosasco15}; see also \cite{LuRudiBorgonovoRosasco20} for recent advances.
It contrasts the approximations \eqref{eq:posterior-mean-approx} and \eqref{eq:posterior-cov-approx} which simply replace the occurrences of $\BFK$ in the posterior distribution of $f$
with $\tilde{\BFK}$.

In particular, it can be shown (see Appendix~\ref{sec:nystrom}) that
\begin{align}
\E[\tilde{f}(\BFx)|\mathcal{D}_n] ={}& \mu(\BFx) + \BFk_m(\BFx)^\intercal (\BFK_{m,m} + \BFK_{m,n}\BFSigma^{-1}\BFK_{n,m})^{-1}  \BFK_{m,n} \BFSigma^{-1}(\bar{\BFy}-\BFmu),   \label{eq:posterior-mean-nystrom}\\
\Cov[\tilde{f}(\BFx),\tilde{f}(\BFx')|\mathcal{D}_n]={}&
K(\BFx,\BFx')-\BFk_m(\BFx)^\intercal (\BFK_{m,m} + \BFK_{m,n}\BFSigma^{-1}\BFK_{n,m})^{-1}  \BFk_m(\BFx'), \label{eq:posterior-cov-nystrom}
\end{align}
both of which can be computed with time complexity $\mathcal{O}(m^2n)$.
We stress, nonetheless, that  \eqref{eq:posterior-mean-nystrom}/\eqref{eq:posterior-cov-nystrom} are not identical to \eqref{eq:posterior-mean-approx}/\eqref{eq:posterior-cov-approx}.

At last, we briefly comment on choosing $m$ and $A$.
Although theoretical analysis may relate the accuracy of the approximation with the asymptotic order of magnitude of $m$ relative to $n$,
in practice $m$ is usually viewed as a tuning parameter.
One may gradually increase the value of $m$, evaluate the resulting accuracy of the approximation, and stop when the marginal improvement falls below some prescribed threshold;
see more discussion in \cite{LuRudiBorgonovoRosasco20}.
Given $m$, $A$ may be determined by random sampling from the entire set of indices $\{1,\ldots,n\}$.


\subsection{Random Features}

Random features represents a large class of algorithms for approximating covariance functions \citep{LiuHuangChenSuykens21}.
We introduce below the original version, called random Fourier features (RFF).
It is proposed in \cite{rahimi2007random}---the seminal work that gives rise to this research direction.

Similar to the Nystr\"om method, RFF also seeks to construct another covariance function $\tilde{K}$ that yields a low-rank approximation to the covariance matrix $\BFK$ to accelerate computation.
However, the approximation that RFF constructs is \emph{data-independent};
that is, it does not depend on the design points $\{\BFx_1,\ldots,\BFx_n\}$.
This is in contrast to the Nystr\"om method, for which the approximations are defined through the design points in the active set.

RFF applies particularly to \emph{stationary} covariance functions, including both the Gaussian and Mat\'ern classes.
If $K$ is stationary covariance function,
then Bochner’s theorem \citep[page~24]{Stein99} asserts that it can be represented as the Fourier transform of a non-negative finite measure:
\begin{equation}\label{eq:Bochner}
K(\BFx, \BFx') = K(\BFzero,\BFzero)\int_{\Real^d} e^{\iu \BFomega^\intercal(\BFx-\BFx')} \mathsf{p}(d\BFomega),
\end{equation}
where $\mathsf{p}(\cdot)$ is a probability measure on $\Real^d$.
For example, if $K$ is the Gaussian covariance function in Equation~\eqref{eq:Gaussian-kernel},
then $\mathsf{p}(\cdot)$ corresponds to the multivariate normal distribution with mean vector $\BFzero$ and covariance matrix $\eta^{-2}\BFI$.
More examples of the $(K,\mathsf{p})$ pair can be found in \cite{LiuHuangChenSuykens21}.

Because $K(\BFx,\BFx')$ is real-valued, we may discard the imaginary part on the right-hand-side of Equation~\eqref{eq:Bochner}. Thus,
\[
K(\BFx, \BFx') = K(\BFzero,\BFzero) \int_{\Real^d} \cos\left( \BFomega^\intercal(\BFx-\BFx')\right) \mathsf{p}(d\BFomega).
\]
Further, if $\BFomega\in\Real^d$ is a random vector having distribution $\mathsf{p}$, then
\begin{align}
K(\BFx,\BFx') ={}& K(\BFzero, \BFzero)\E_{\BFomega}\left[\cos\left( \BFomega^\intercal(\BFx-\BFx')\right)\right]  \nonumber \\
={}& K(\BFzero, \BFzero)\E_{\BFomega, b}\left[\sqrt{2}\cos\left( \BFomega^\intercal \BFx + b\right) \sqrt{2}\cos\left( \BFomega^\intercal \BFx' + b\right)\right], \label{eq:cos-prod}
\end{align}
where $b$ is an independent random variable uniformly distributed on $(0,2\pi)$. The proof of Equation~\eqref{eq:cos-prod} is provided in Appendix~\ref{sec:random-feat}.

We then apply the standard Monte Carlo approximation:
\begin{align*}
K(\BFx,\BFx')  \approx{}& K(\BFzero,\BFzero)\cdot\frac{1}{m}\sum_{t=1}^m \sqrt{2}\cos\left( \BFomega_t^\intercal \BFx + b_t\right) \sqrt{2}\cos\left( \BFomega_t^\intercal \BFx' + b_t\right) \\
={}&
\begin{pmatrix}
\sqrt{\frac{2K(\BFzero,\BFzero)}{m}} \cos\left( \BFomega_1^\intercal \BFx + b_1\right) &
\cdots &
\sqrt{\frac{2K(\BFzero,\BFzero)}{m}} \cos\left( \BFomega_m^\intercal \BFx + b_m\right)
\end{pmatrix}
\begin{pmatrix}
\sqrt{\frac{2K(\BFzero,\BFzero)}{m}} \cos\left( \BFomega_1^\intercal \BFx' + b_1\right) \\
\vdots \\
\sqrt{\frac{2K(\BFzero,\BFzero)}{m}} \cos\left( \BFomega_m^\intercal \BFx' + b_m\right)
\end{pmatrix} \\
\coloneqq  {}&  \BFphi_m(\BFx)^\intercal \BFphi_m(\BFx') \coloneqq \tilde{K}(\BFx,\BFx'),
\end{align*}
where $\{\omega_t:t=1,\ldots,m\}$ are independent samples drawn from $\mathsf{p}$, and
$\{b_t:t=1,\ldots,m\}$ are independent samples drawn from $\mathsf{Uniform}(0,2\pi)$.
It is easy to show that $\tilde{K}$ is a covariance function.

Note that, in light of the discussion in Section~\ref{sec:connect},
$\BFphi_m(\BFx)$ can be view as a vector of basis functions---also known as \emph{features} in the machine learning literature---and they are constructed via random sampling, hence the method's name ``random features.''

Let $\tilde{\BFK} \coloneqq (\tilde{K}(\BFx_i,\BFx_{i'}))_{i,i'=1}^n$.
Then,
\[
\tilde{\BFK} = \begin{pmatrix}
\BFphi_m(\BFx_1)^\intercal \\
\vdots \\
\BFphi_m(\BFx_n)^\intercal
\end{pmatrix}
\begin{pmatrix}
\BFphi_m(\BFx_1) &
\cdots &
\BFphi_m(\BFx_n)
\end{pmatrix}
\coloneqq
\BFPhi_m \BFPhi_m ^\intercal.
\]
Because $\BFPhi_m$ is a $n$-by-$m$ matrix,
$\tilde{\BFK}$ is a low-rank approximation of the covariance matrix $\BFK$, provided that $m<n$.
In addition,
let $\tilde{f}\sim\mathsf{GP}(\mu, \tilde{K})$. We prove in Appendix~\ref{sec:random-feat} that
\begin{align}
\E[\tilde{f}(\BFx)|\mathcal{D}_n] ={}& \mu(\BFx) + \BFphi_m(\BFx)^\intercal (\BFI + \BFPhi_m^\intercal \BFSigma^{-1}\BFPhi_m)^{-1}  \BFPhi_m^\intercal \BFSigma^{-1}(\bar{\BFy}-\BFmu),  \label{eq:posterior-mean-RFF}\\
\Cov[\tilde{f}(\BFx),\tilde{f}(\BFx')|\mathcal{D}_n]={}&
K(\BFx,\BFx')-\BFphi_m(\BFx)^\intercal (\BFI + \BFPhi_m^\intercal \BFSigma^{-1}\BFPhi_m)^{-1}  \BFphi_m(\BFx'). \label{eq:posterior-cov-RFF}
\end{align}
Similar to the Nystr\"om approximations \eqref{eq:posterior-mean-nystrom}--\eqref{eq:posterior-cov-nystrom},
the time complexity for computing \eqref{eq:posterior-mean-RFF}--\eqref{eq:posterior-cov-RFF} is also $\mathcal{O}(m^2n)$.

\section{Concluding Remarks}\label{sec:conclusions}
We have introduced several common surrogates---low-order polynomials, linear basis function models, and Gaussian processes---along with  simple techniques to enhance their prediction capability with little additional computational overhead.
With the use of these surrogates, we have presented a number of approaches to solving SO problems with continuous decision variables, some of which (RSM, STRONG, SPAS) are locally convergent while others (KG, GP-UCB, GPS) globally convergent. Moreover, we have discussed two widely popular methods---the Nystr\"om method and random features---for dealing with computational challenges associated with Gaussian processes in the presence of large datasets.

Looking forward, we believe the following research directions are potentially fruitful and of high impact.
First, an ideal SO algorithm in our mind would be able to quickly identify most local optima and then select the best among them.
In contrast, it is often observed in practice that globally convergent methods such as those introduced in this tutorial tend to ``over-explore''---direct sampling efforts away from the proximity of a local optimum without realizing it.
Lack of gradient information to provide local curvature of the response surface, among others, is an important reason.
It is of great interest to develop algorithms that
integrate local search and global search to ensure fast convergence to global optima.

Second, another plausible approach to accelerating convergence to global optima is to incorporate structural information---such as monotonicity, convexity, and level of differentiability---provided that one can safely impose such assumptions on the response surface.
Surrogates that process such properties do exist \citep{LimGlynn12,WangBerger16,SalemiStaumNelson19}.
But how to leverage them to develop fast SO algorithms is yet to be fully investigated.
See \cite{zhang2020discrete} and \cite{zhang2020stochastic} for recent developments in this regard.

Third,
most theoretical analyses of SO algorithms in the literature focus on proving convergence to a local/global optimum as the computational budget grows.
Such asymptotic analysis can hardly explain the algorithms' performance in finite time,
nor can it provide accurate guidance for performance tuning in practice.
Little is known about their rates of convergence.
Notable exceptions include \cite{Bull11} and \cite{BouttierGavra19}.
The former establishes the rate of convergence for Efficient Global Optimization algorithms \citep{JonesSchonlauWelch98},
while the latter for Simulated Annealing algorithms \citep{GelfandMitter89}.
In general, deeper theoretical understanding of SO algorithms is strongly needed to fill the gap,
which may shed light on critical attributes required  for improving algorithm efficiency.

Fourth, in Section~\ref{sec:STRONG} we introduced the STRONG algorithm that integrates low-order polynomial surrogates with trust-region methods. Even though there are a number of papers taking this approach, we think it is under-studied and has a potential to solve large-scale SO problems. To build a full quadratic model, one needs at least $d(d+1)/2+1$ design points and each design point may need more than one replication of the simulation experiments in order to take into consideration the simulation noise. If one needs to conduct these many experiments in each iteration, it may be too costly. One way that may solve the problem is to use $L_1$-regularization, a.k.a. LASSO, in the regression to select the most important parameters of the quadratic model with far fewer experiments \citep{Tibshirani96}. Notice that these selected parameters include an important part of the ascent and curvature information and one may use this partial information to guide the optimization process. Even though we think the idea of adding $L_1$-regularization into the trust-region framework has great potential for SO, it is quite challenging to design efficient algorithms and to analyze their asymptotic properties, e.g., convergence and rate of convergence.

Fifth, running simulation experiments is often time-consuming. But different simulation experiments are typically independent and they can be run on different processors. Therefore, it is natural to think how to design SO algorithms that work well in parallel computing environments. Recently, \cite{LuoHongNelsonWu15}, \cite{NiCiocanHendersonHunter17}, \cite{ZhongHong21}, and others have developed parallel algorithms for R\&S problems.
\cite{WuFrazier16} also develop parallel knowledge gradient algorithms for Bayesian optimization. However, we believe there are still many opportunities in developing efficient surrogate-based parallel SO algorithms.

Last but not the least, in the present era of big data an emerging decision-making paradigm that has gained great popularity in recent years is optimization with covariates.
That is, the optimal decision is no longer constant, but varies as a function of the covariates that represent the additional contextual information available at the moment of making a decision; see, e.g., \cite{BanRudin19}, \cite{BertsimasKallus20}, and \cite{BertsimasKoduri21}.
However, these articles focus on settings where
the objective functions are analytically tractable.
\cite{ShenHongZhang21} address the problem of R\&S with covariates, assuming the response surface  for each alternative is a linear function in the covariates.
We expect much more to be explored in the direction of SO with covariates for years to come.

%
%
%
\begin{APPENDICES}
\section{Proofs Related to the Nystr\"om Method}\label{sec:nystrom}
Let
$\tilde{K}(\BFx,\BFx') \coloneqq  \BFk_m(\BFx)^\intercal \BFK_{m,m}^{-1} \BFk_m(\BFx')$,
$\tilde{\BFk}(\BFx)\coloneqq (\tilde{K}(\BFx,\BFx_1),\ldots,\tilde{K}(\BFx,\BFx_n))^\intercal$, and
$\BFQ \coloneqq\BFK_{m,m} + \BFK_{m,n}\BFSigma^{-1}\BFK_{n,m}$.
Then,
\begin{align}
\tilde{\BFk}(\BFx)^\intercal(\tilde{\BFK}+\BFSigma)^{-1}
={}& \BFk_m(\BFx)^\intercal \BFK_{m,m}^{-1} \BFK_{m,n} (\BFK_{n,m}\BFK_{m,m}^{-1}\BFK_{m,n} + \BFSigma)^{-1} \nonumber \\
={}& \BFk_m(\BFx)^\intercal \BFK_{m,m}^{-1} \BFK_{m,n} [\BFSigma^{-1} - \BFSigma^{-1} \BFK_{n,m}(\BFK_{m,m} + \BFK_{m,n}\BFSigma^{-1}\BFK_{n,m})^{-1}\BFK_{m,n}\BFSigma^{-1} ]\label{eq:applying-woodbury}\\
={}& \BFk_m(\BFx)^\intercal \BFK_{m,m}^{-1}  (\BFK_{m,n}\BFSigma^{-1} - \BFK_{m,n}\BFSigma^{-1} \BFK_{n,m}\BFQ^{-1}\BFK_{m,n}\BFSigma^{-1} ) \nonumber \\
={}& \BFk_m(\BFx)^\intercal \BFK_{m,m}^{-1} (\BFI - \BFK_{m,n}\BFSigma^{-1}\BFK_{n,m}\BFQ^{-1})\BFK_{m,n}\BFSigma^{-1} \nonumber\\
={}& \BFk_m(\BFx)^\intercal \BFK_{m,m}^{-1} (\BFQ - \BFK_{m,n}\BFSigma^{-1}\BFK_{n,m})\BFQ^{-1}\BFK_{m,n}\BFSigma^{-1}\nonumber\\
={}& \BFk_m(\BFx)^\intercal\BFQ^{-1}\BFK_{m,n}\BFSigma^{-1}, \label{eq:woodbury-2}
\end{align}
where Equation~\eqref{eq:applying-woodbury} follows from the Woodbury matrix identity.

Applying Equation~\eqref{eq:posterior-mean} to $\tilde{f}\sim \mathsf{GP}(\mu, \tilde{K})$, we have
\begin{align*}
\E[\tilde{f}(\BFx)|\mathcal{D}_n] ={}& \mu(\BFx)   + \tilde{\BFk}(\BFx)^\intercal(\tilde{\BFK}+\BFSigma)^{-1} (\bar{\BFy} - \BFmu)  \\
={}&  \mu(\BFx)   + \BFk_m(\BFx)^\intercal\BFQ^{-1}\BFK_{m,n}\BFSigma^{-1} (\bar{\BFy} - \BFmu) \\
={}&  \mu(\BFx)   + \BFk_m(\BFx)^\intercal (\BFK_{m,m} + \BFK_{m,n}\BFSigma^{-1}\BFK_{n,m})^{-1}  \BFK_{m,n} \BFSigma^{-1} (\bar{\BFy} - \BFmu),
\end{align*}
where the second equation follows from \eqref{eq:woodbury-2}.
This completes the proof of Equation~\eqref{eq:posterior-mean-nystrom}.

Also by Equation~\eqref{eq:woodbury-2},
\begin{align*}
\tilde{\BFk}(\BFx)^\intercal  (\tilde{\BFK}+\BFSigma)^{-1} \tilde{\BFk}(\BFx')
={}&
\BFk_m(\BFx)^\intercal\BFQ^{-1}\BFK_{m,n}\BFSigma^{-1} \tilde{\BFk}(\BFx') \\
={}& \BFk_m(\BFx)^\intercal\BFQ^{-1}\BFK_{m,n}\BFSigma^{-1} \BFK_{n,m}\BFK_{m,m}^{-1}\BFk_m(\BFx').
\end{align*}
It follows that, after applying Equation~\eqref{eq:posterior-cov} to $\tilde{f}\sim \mathsf{GP}(\mu, \tilde{K})$,
\begin{align*}
K(\BFx,\BFx')-\Cov[\tilde{f}(\BFx),\tilde{f}(\BFx')|\mathcal{D}_n] ={}& \tilde{\BFk}(\BFx)^\intercal  (\tilde{\BFK}+\BFSigma)^{-1} \tilde{\BFk}(\BFx') \\
={}& \BFk_m(\BFx)^\intercal \BFK_{m,m}^{-1} \BFk_m(\BFx)
- \BFk_m(\BFx)^\intercal\BFQ^{-1}\BFK_{m,n}\BFSigma^{-1}  \BFK_{n,m}\BFK_{m,m}^{-1}\BFk_m(\BFx')
\\
={}& \BFk_m(\BFx)^\intercal (\BFI - \BFQ^{-1}\BFK_{m,n}\BFSigma^{-1}  \BFK_{n,m} )\BFK_{m,m}^{-1} \BFk_m(\BFx')
\\
={}& \BFk_m(\BFx)^\intercal \BFQ^{-1}(\BFQ - \BFK_{m,n}\BFSigma^{-1}  \BFK_{n,m})\BFK_{m,m}^{-1} \BFk_m(\BFx') \\
={}& \BFk_m(\BFx)^\intercal \BFQ^{-1} \BFk_m(\BFx') \\
={}& \BFk_m(\BFx)^\intercal (\BFK_{m,m} + \BFK_{m,n}\BFSigma^{-1}\BFK_{n,m})^{-1}  \BFk_m(\BFx'),
\end{align*}
proving Equation~\eqref{eq:posterior-cov-nystrom}.

\section{Proofs Related to Random Features}\label{sec:random-feat}
Suppose $b\sim \mathsf{Uniform}(0,2\pi)$.
Then, for any $a\in\Real$,
\begin{align*}
\E_b[\cos(a+2b)] =  \int_0^{2\pi} \frac{\cos(a + 2b)}{2\pi}\, db = \frac{1}{2\pi}\sin(a + 2b) \Big|_{0}^{2\pi } = \frac{1}{2\pi} [\sin(a+4\pi) -\sin(a) ] = 0.
\end{align*}
Thus,
\[
\E_{\BFomega, b} \bigl[\cos\bigl(\BFomega^\intercal(\BFx+\BFx') + 2b) ] =
\E_{\BFomega}\bigl[\E_b\bigl[\cos\bigl(\BFomega^\intercal(\BFx+\BFx'\bigr) + 2b)\big|\BFomega\bigr]\bigr] = 0.
\]
It follows that
\begin{align*}
\E_{\BFomega, b}\bigl[\cos\bigl(\BFomega^\intercal(\BFx-\BFx')\bigr)\bigr] = {}&
\E_{\BFomega, b}\bigl[\cos\bigl(\BFomega^\intercal(\BFx+\BFx') + 2b\bigr) \bigr] + \E_{\BFomega, b}\bigl[\cos\bigl(\BFomega^\intercal(\BFx-\BFx')\bigr)\bigr] \\
={}&  \E_{\BFomega, b}\bigl[\cos\bigl((\BFomega^\intercal\BFx + 2b) + (\BFomega^\intercal\BFx' + 2b) \bigr) \bigr] +  \E_{\BFomega, b}\bigl[\cos\bigl((\BFomega^\intercal\BFx + 2b) - (\BFomega^\intercal\BFx' + 2b) \bigr) \bigr]
\\
={}& \E_{\BFomega, b} \bigl[2\cos( \BFomega^\intercal \BFx + b) \cos( \BFomega^\intercal \BFx' + b)\bigr],
\end{align*}
proving Equation~\eqref{eq:cos-prod}.

We now prove Equations~\eqref{eq:posterior-mean-RFF} and \eqref{eq:posterior-cov-RFF},
following the strategy in Appendix~\ref{sec:nystrom}.
Let $\tilde{K}(\BFx,\BFx') \coloneqq \BFphi_m(\BFx)^\intercal \BFphi_m(\BFx')  $ and $\BFQ\coloneqq \BFI + \BFPhi_m^\intercal\BFSigma^{-1}\BFPhi_m$.
The key is to derive $\tilde{\BFk}(\BFx)^\intercal(\tilde{\BFK}+\BFSigma)^{-1}$ as follows.
\begin{align*}
\tilde{\BFk}(\BFx)^\intercal(\tilde{\BFK}+\BFSigma)^{-1}
={}& \BFphi_m(\BFx)^\intercal  \BFPhi_m^\intercal (\BFPhi_m \BFPhi_m^\intercal + \BFSigma)^{-1}  \\
={}& \BFphi_m(\BFx)^\intercal  \BFPhi_m^\intercal [\BFSigma^{-1} - \BFSigma^{-1} \BFPhi_m (\BFI + \BFPhi_m^\intercal\BFSigma^{-1}\BFPhi_m )^{-1}\BFPhi_m^\intercal\BFSigma^{-1} ]\label{eq:applying-woodbury}\\
={}& \BFphi_m(\BFx)^\intercal   (\BFPhi_m^\intercal\BFSigma^{-1} - \BFPhi_m^\intercal\BFSigma^{-1} \BFPhi_m \BFQ^{-1}\BFPhi_m^\intercal\BFSigma^{-1} )  \\
={}& \BFphi_m(\BFx)^\intercal  (\BFI - \BFPhi_m^\intercal\BFSigma^{-1}\BFPhi_m \BFQ^{-1})\BFPhi_m^\intercal\BFSigma^{-1} \\
={}& \BFphi_m(\BFx)^\intercal  (\BFQ - \BFPhi_m^\intercal\BFSigma^{-1}\BFPhi_m )\BFQ^{-1}\BFPhi_m^\intercal\BFSigma^{-1}\\
={}& \BFphi_m(\BFx)^\intercal\BFQ^{-1}\BFPhi_m^\intercal\BFSigma^{-1}.
\end{align*}
The remaining calculations are essentially the same as those in Appendix~\ref{sec:nystrom}, so we omit them.

\end{APPENDICES}


\ACKNOWLEDGMENT{L. Jeff Hong is supported in part by the Natural Science Foundation of China (Project No. 72091211 and Project No. 71991473). Xiaowei Zhang is supported in part by the Hong Kong Research Grant Council (Project No. 16211417 and Project No. 17201520).}


\bibliographystyle{informs2014} 
\bibliography{ref.bib} 






\end{document}